\documentclass[11pt,leqno]{article}
\usepackage[centertags]{amsmath}
\usepackage{amsfonts}
\usepackage{amsmath}
\usepackage{amssymb}
\usepackage{latexsym}
\usepackage{amsthm}
\usepackage{mathrsfs}
\usepackage{newlfont}
\usepackage[latin1]{inputenc}

\newtheorem{teor}{Theorem}[section]

\newtheorem{lema}{Lemma}[section]

\newcommand\CC{{\mathbb C}}

\newcommand\RR{{\mathbb R}}

\newcommand\NN{{\mathbb N}}

\newcommand\X{{\mathfrak X}}
\newcommand\tr{\operatorname{Tr}}

\newcommand{\D}{\displaystyle}

\numberwithin{equation}{section} 

\title{Second order differential operators having several families
of orthogonal matrix polynomials as eigenfunctions\footnote{The
  work of the  authors is  partially supported by  D.G.E.S,
  ref. BFM2006-13000-C03-01, FQM-262, P06-FQM-01735, FQM-481 (\textit {Junta de
Andaluc\'{\i}a)}. \newline \textit{2000 Mathematics Subject
Classification}. Primary 42C05. \newline \textit{Key words and
phrases}. Orthogonal matrix polynomials, Differential equations.}}
\author{Antonio J. Dur\'an and Manuel D. de la Iglesia$^{\dagger}$ \\
   \footnotesize $\dagger$  \footnotesize
    \  Departamento de An\'{a}lisis Matem\'{a}tico.
   Universidad de Sevilla \\
   \footnotesize Apdo (P. O. BOX) 1160. 41080 Sevilla. Spain.
   duran@us.es, mdi29@us.es \\}
\date{}

\begin{document}

\maketitle

\begin{abstract}
The aim of this paper is to bring into the picture a new phenomenon
in the theory of orthogonal matrix polynomials satisfying second
order differential equations. The last few years have witnessed some
examples of a (fixed) family of orthogonal matrix polynomials whose
elements are common eigenfunctions of several linearly independent
second order differential operators. We show that the dual situation
is also possible: there are examples of one parametric families of
monic matrix polynomials, each family orthogonal with respect to a
different weight matrix, whose elements are eigenfunctions of a
common second order differential operator.

These examples are constructed by adding a discrete mass to a weight matrix at a certain point. In this article it is described how to choose a point $t_0$, a discrete mass $M(t_0)$ and the weight matrix $W$ so that the new weight matrix $W+\delta_{t_0}M(t_0)$ inherits some of the
symmetric second order differential operators associated with $W$.
It is well known that this situation is not possible for the
classical scalar families of Hermite, Laguerre and Jacobi.

For some of these examples we characterize the convex cone of weight
matrices for which the differential operator is symmetric.
\end{abstract}

\section{Introduction} The theory of matrix valued orthogonal polynomials starts with two
papers by M. G. Kre$\breve{\mbox{{\i}}}$n in 1949, see \cite{K1,
K2}. A sequence of orthonormal matrix polynomials $(P_n)_n$ can be
characterized as solutions of the difference equation
\begin{align*}
tP_n(t)=A_{n+1}P_{n+1}(t)+B_nP_n(t)+A_n^*P_{n-1}(t),\quad
n=0,1,\ldots,
\end{align*}
where $A_n$ and $B_n$ are $N\times N$ nonsingular and
Hermitian matrices, respectively, and initial conditions $P_{-1}=0$ and $P_0$
nonsingular. Each family $(P_n)_n$ goes along with a weight matrix
$W$ and satisfies $\int P_ndWP_m^*=\delta _{n,m}I$.

More than 50 years later the first examples of orthogonal matrix polynomials $(P_n)_n$ satisfying
second order differential equations of the form
\begin{align}\label{difequ}
P''_n(t)F_2(t)+P'_n(t)F_1(t)+P_n(t)F_0=\Gamma _n P_n(t), \quad
n=0,1,\ldots,
\end{align}
were produced. Here $F_2$, $F_1$ and $F_0$ are matrix polynomials (which do not
depend on $n$) of degrees less than or equal to $2$, $1$ and $0$,
respectively (see \cite{DG1,G1,GPT}). Two main methods have been
developed in the last five years to produce such examples: solving
an appropriate set of differential equations (see \cite{D2, DG1,
DG5, DdI}) or coming from the study of matrix valued spherical
functions (see \cite{GPT1, GPT6, PT1}). In the case of matrix orthogonality, these families of orthogonal
matrix polynomials are likely to play the role of the classical families of Hermite, Laguerre and Jacobi in the case of scalar orthogonality. The complexity of the matrix world however proved to be much richer compared to the scalar case (see, for instance, the papers cited above).

We point out here that, if the eigenvalues
$\Gamma _n$ are Hermitian, then the second order differential equation
(\ref{difequ}) for the orthonormal polynomials $(P_n)_n$ is
equivalent to the symmetry of the second order differential operator
\begin{align}\label{difope}
D=\partial ^2F_2(t)+\partial^1F_1(t)+\partial^0F_0,\quad \partial
=\frac{d}{dt},
\end{align}
with respect to the weight matrix $W$ (see \cite{D1}). The symmetry of $D$ with
respect to $W$ is defined as
\begin{align}\label{rev1}
\int (P D) dW Q^*=\int P dW (Q D)^*,
\end{align}
for any matrix polynomials $P$ and $Q$. Here and in the rest of the
paper we will follow the notation in \cite{GT} for right-hand side
differential operators. In particular, if $P$ is a matrix polynomial
and $D$ a differential operator as (\ref{difope}), by $PD$ we mean
$$
PD=P''(t)F_2(t)+P'(t)F_1(t)+P(t)F_0.
$$

As more families of orthogonal matrix polynomials satisfying second
order differential equations become available many new interesting phenomena are being discovered, which are
absent in the well known scalar theory.

One of such phenomena is the fact that the elements of a family of
orthogonal matrix polynomials $(P_n)_n$ can be common eigenfunctions
of several linearly independent second order differential operators
(while in the scalar case the symmetric second order differential operator is unique up to multiplicative and additive constants). The first illustrations of this phenomenon have been recently found by F. A. Gr\"unbaum and M. M.
Castro \cite{CG1} and other authors contributed more later. Some of the examples arise from group representation theory. For instance, \cite{GPT1} discusses two second order
differential operators acting on matrix spherical functions which
were later put in the framework of orthogonal polynomials in \cite{GPT3,
GPT6, PT1, PR}; see also \cite{GdI}. Other examples were found by integrating an appropriate set of differential equations (see \cite{D2, DdI, DL}). See also \cite{CG2}, where the authors take up
the issue of existence of orthogonal matrix polynomials which are
common eigenfunctions of differential operators of order one.

As a consequence of this phenomenon the algebra of differential
operators associated with a fixed weight matrix $W$ is receiving a lot of attention. This algebra $\mathcal {D}(W)$ is
defined as follows: given a fixed sequence of orthogonal polynomials
$(P_n)_n$ with respect to  $W$ (the monic sequence, for instance),
$\mathcal {D}(W)$ is formed by all differential operators
\begin{equation}\label{dop1}
D=\sum_{i=0}^k \partial^i F_i(t),\quad \partial=\frac{d}{dt},\quad k\geq0,
\end{equation}
where $F_i(t), i=0,\ldots,k$ are matrix polynomials of $\deg
(F_i)\le i$, for which $P_nD=\Gamma_n P_n$, $n\ge 0$. For the
classical families of the scalar case every differential operator
having one of the families as eigenfunctions has to be a polynomial in the
corresponding symmetric second order differential operator (see
\cite{M}). Hence the associated algebra is isomorphic to $\CC[t]$.
The examples of weight matrices having several linearly independent
symmetric second order differential operators show that in the
matrix case the problem of characterizing the algebra $\mathcal
{D}(W)$ is going to be a rather more difficult problem. In
\cite{CG1, GdI, DdI} some conjectures have been made about the structure of the
algebra $\mathcal D(W)$ for some concrete examples.
Based on computational evidence these conjectures show that we
can expect a big variety of situations in the matrix case. To the best of our knowledge, only one of those conjectures
has been proved for the weight matrix (\ref{51dg}) (see
\cite{CG1}, Sect. 6, for the conjecture and \cite{T} for the proof).

The purpose of this paper is to show what one can call the dual
situation to that described in the previous paragraph. For a fixed
differential operator $D$ of the form (\ref{dop1}), we define a
set of weight matrices
\begin{align}\label{yd}
\Upsilon (D)=\{ W: \mbox{$D$ is symmetric with respect to $W$}\},
\end{align}
where the symmetry of $D$ is defined again by (\ref{rev1}).

Note that if $\Upsilon (D)\not =\emptyset$ then
it is a convex cone:
  if $W_1, W_2\in \Upsilon (D)$ and $\gamma,  \zeta \ge 0$
  (one of them non null), then $\gamma W_1+\zeta W_2\in \Upsilon (D)$.

The weight matrices $W$ going along with a symmetric second order
differential operator $D$ mentioned at the beginning of this paper
provide examples where $\Upsilon (D)\not =\emptyset$. In these
examples $\Upsilon (D)$ contains at least a half line: $\gamma W$,
$\gamma >0$. In this paper we show the first examples of operators
$D$ for which $\Upsilon (D)$ is a two dimensional convex cone. That
is, we show examples of a fixed second order differential operator
$D$ as (\ref{difope}) for which there exist two weight matrices
$W_1$ and $W_2$, $W_1\not =\alpha W_2$ for any $\alpha >0$, such
that $D$ is symmetric with respect to any of the weight matrices
$\gamma W_1+\zeta W_2$, $\gamma,  \zeta \ge 0 $. That means, in
particular, that the corresponding monic matrix polynomials
$(P_{n,\zeta /\gamma })_n$ orthogonal with respect to  $\gamma
W_1+\zeta W_2$ (they only depend on $W_1$, $W_2$ and the ratio
$\zeta/\gamma$) are eigenfunctions of $D$
$$
P_{n,\zeta /\gamma}D=\Gamma_nP_{n,\zeta /\gamma },\quad n=0, 1,
\ldots,\quad \gamma >0,\quad \zeta \geq0,
$$
where $D$ and $\Gamma_n$ do not depend on $\gamma, \zeta$.

We give a simple but fruitful method to find such examples (Section
\ref{sec1}) and show a collection of instructive examples (Section
\ref{sec2}).

Our method itself is a surprise if one compares it with the
situation in the scalar case. We first take a weight matrix $W$
which has several linearly independent symmetric second order
differential operators. And then we add a Dirac distribution
$\delta _{t_0}M(t_0)$ to $W$, where the real number $t_0$ and the mass
$M(t_0)$ (a Hermitian positive semidefinite matrix) are carefully
chosen. We show in Section 2 that for a fixed $t_0$ and under certain mild conditions, we can produce a positive semidefinite matrix $M(t_0)$ and a second order differential operator $D$ symmetric with
respect to $W$, such that $D$ is also symmetric with respect to any
weight matrix of the form $\gamma W+\zeta\delta _{t_0}M(t_0)$,
$\gamma>0, \zeta \ge 0$. In Section 3 we illustrate with some examples that the choice of the point $t_0$ where the Dirac distribution is located depends more on the matrices $F_2(t_0)$,
$F_1(t_0)$ and $F_0$ than on the support of $W$. We also
characterize for some of these examples the convex cone (\ref{yd}).

The situation is quite different from the scalar case.
When a mass point is added to any of the classical weights of Hermite,
Laguerre and Jacobi, the existence of a symmetric second order
differential operator automatically disappears. Only when $t_0$ is
taken at the endpoints of the support one eventually gets the
symmetry of a fourth (or even larger) order differential operator
which is not symmetric with respect to the original weight. This
arises for the particular cases of the Laguerre weight $e^{-t}$ in
$(0,+\infty )$, for the Legendre weight $1$ in $(-1,1)$ and for the
special case of the Jacobi weight $(1-t)^\alpha$ in $(0,1)$, raising
the so called Laguerre type weight $e^{-t}+M\delta _0$, Legendre
type weight $1+M(\delta _{-1}+\delta_1)$ and Jacobi type weight
$(1-t)^\alpha+M\delta _0$, respectively (see \cite{LK1}, or
\cite{GH} where these and some other examples are obtained by
applying the Darboux process).

We would like to conclude this Introduction by
displaying one of our examples. Consider the weight matrix
\begin{equation}\label{51dg}
W_a(t)=e^{-t^2}\begin{pmatrix} 1+a^2t^2 & at \\
                                      at & 1 \\
                                    \end{pmatrix}
,\quad t\in\mathbb{R},\quad a\in\mathbb{R}\setminus \{0\}.
\end{equation}
The linear space of differential operators of order at most two
having the orthogonal polynomials with respect to $W_a$ as
eigenfunctions has dimension five. A basis is formed by the identity
and four linearly independent operators of order two (see Section 6
of \cite{CG1}). We show in Section 3.1 that the weight matrices
$W_{a,\gamma,\zeta}=\gamma W_a+\zeta\delta_0\begin{pmatrix}
      1& 1 \\
      1 &1 \\
    \end{pmatrix}$, $\gamma>0, \zeta \geq0$, share the following symmetric second order differential
    operator:
\begin{equation*}
D_a=\partial^2\begin{pmatrix}
      1-at & -1+a^2t^2 \\
      -1 & 1+at \\
    \end{pmatrix}+\partial^1\begin{pmatrix}
      -2a-2t & 2a+2(2+a^2)t \\
      0& -2t \\
    \end{pmatrix}+\partial^0\begin{pmatrix}
      -1 & 2\D\frac{2+a^2}{a^2} \\
      \D\frac{4}{a^2} & 1 \\
    \end{pmatrix}.
\end{equation*}
Moreover, we prove that
$$
\Upsilon (D_a)=\left\{ \gamma W_a+\zeta \delta_0\begin{pmatrix}
      1& 1 \\
      1 &1 \\
    \end{pmatrix}: \gamma >0, \zeta \ge 0\right\}.
$$
As we pointed out above, this means that each  monic family
$(P_{n,a,\zeta/\gamma})_n$, $\gamma>0, \zeta\geq0$, orthogonal with
respect to $W_{a,\gamma,\zeta}$  satisfies the same second order
differential equation, namely
\begin{equation*}
P_{n,a,\zeta/\gamma}D_a=\Gamma_{n,a}P_{n,a,\zeta/\gamma},\quad n=0,
1, \ldots,
\end{equation*}
where
$$
\Gamma_{n,a}=\begin{pmatrix}
      -(2n+1)& \D\frac{(2+na^2)(2+(n+1)a^2)}{a^2} \\
      \D\frac{4}{a^2} & -2n+1
    \end{pmatrix}.
$$
Notice that neither $D_a$ nor $\Gamma_{n,a}$ depend on
$\gamma,\zeta$.

We will also see in Section \ref{ssH} that the point $t_0$ where
the discrete mass is added can be located at any real number. There
always exists a Hermitian positive semidefinite matrix $M(t_0)$ such
that the weight matrices $\gamma W_a+\zeta \delta_{t_0}M(t_0)$,
$\gamma>0,\zeta\geq0$, share a common symmetric second order
differential operator (see (\ref{dopherm})).

\bigskip

\section{The main result}\label{sec1}
Presented in this section is a set of constraints to guarantee
the symmetry of a differential operator of any order with respect to
a weight matrix modified by adding  a Dirac distribution at an
arbitrary point. Before that we need some definitions and previous
results.

We call an $N\times N$ matrix of measures  $W$ (supported in the
real line) a  weight matrix if it satisfies the following: the numerical
matrix $W(\Omega )$ is positive semidefinite for any Borel set
$\Omega $; the integral $\mu_n=\int  t^ndW(t)$ (called the $n$-th
moment of $W$) exists and is finite for any $n\in \NN$; and $\int
P(t)\,dW(t)\,P^*(t)$ is nonsingular if the leading coefficient of
the matrix polynomial $P$ is nonsingular.

A Hermitian sesquilinear form in the linear space of matrix
polynomials can be associated with a weight matrix $W$:
\begin{equation*}
    \langle P,Q\rangle=\int P(t)\,dW(t)\,Q^*(t),
\end{equation*}
where $Q^*(t)$ denotes the conjugate transpose of $Q(t)$.

We can now produce a sequence of orthogonal matrix polynomials
$(P_n)_n$ with deg $P_n=n$ and nonsingular leading coefficient such
that $\langle P_n,P_m\rangle =\Delta_n\delta _{n,m}$ with
$\Delta_n$ positive definite. If $\Delta_n=I$ we say that the
sequence $(P_n)_n$ is orthonormal (for a much more complete
introduction to matrix orthogonality, see \cite{DG2} and references
therein).

\bigskip

If one is considering possible applications of orthogonal matrix
polynomials, it is natural to concentrate on the cases where some
additional property holds. For instance, in \cite{D1} one of the
authors raised a problem of characterizing weight matrices whose
orthonormal matrix polynomials are common eigenfunctions of some
symmetric right-hand side second order differential operator $D$ as
(\ref{difope}) with Hermitian eigenvalues. We say that a
differential operator $D$ is symmetric with respect to a weight
matrix $W$ if
$$
\langle PD,Q\rangle =\langle P,QD\rangle
$$
for any pair of matrix polynomials $P$ and $Q$. Recall that we
are following the notation in \cite{GT} for right-hand side
differential operators. We already mentioned in the Introduction
that in the last few years a large class of families of weight
matrices $W$ has been found having symmetric second order
differential operators as (\ref{difope}).

The condition of symmetry for the pair made up of a weight matrix
$W$ and a  differential operator $D_k$ of order $k$ can be
established in terms of a set of difference and differential
equations relating $W$ and the coefficients of $D_k$. Indeed, if we
write the right-hand side differential operator $D_k$ of order $k$
as
\begin{equation}\label{dop}
D_k=\sum_{i=0}^k \partial^i F_i(t),\quad \partial=\frac{d}{dt},
\end{equation}
where $F_i(t), i=0,\ldots,k$ are matrix polynomials of degree less
than or equal to $i$,
$$
F_i(t)=\sum_{j=0}^i t^jF_{j}^i ,\;\;\; F_{j}^i\in
\mathbb{C}^{N\times N},
$$
and denote by $\mu_n$, $n=0,1,\ldots$, the moments of the weight
matrix $W$,  then we have the following:

\begin{teor}\label{teorsymm}
For a weight matrix $W$ the following two conditions are equivalent:
\begin{enumerate}
  \item The operator $D_k$ is symmetric with respect to $W$.
  \item For $n\geq l$, the following $k+1$ sets of moment equations
  hold
  \begin{equation}\label{moment}
\sum_{i=0}^{k-l}
\begin{pmatrix}
k-i \\
l \\
\end{pmatrix}
(n-l)_{k-l-i}B_n^{k-i}=(-1)^l(B_n^l)^*,\quad l=0,\ldots,k,
\end{equation}
where
$$
B_n^l=\sum_{i=0}^l F_{l-i}^l\mu_{n-i},\quad l=0,\ldots,k.
$$
\end{enumerate}
Moreover, suppose the weight matrix $W=W(t)dt$ has a smooth density
$W(t)$ with respect to the Lebesgue measure which satisfies the
boundary conditions that
    \begin{equation}\label{bound}
    \sum_{i=0}^{p-1}(-1)^{k-i+p-1}
    \begin{pmatrix}
    k-i \\
    l \\
    \end{pmatrix}
    \big(F_{k-i}\cdot W\big)^{(p-1-i)},\quad
    p=1,\ldots,k, \quad l=0,\ldots,k-p,
    \end{equation}
    should have vanishing limits at each of the endpoints of the support
    of $W$, and the following $k+1$
    matrix differential equations hold
    \begin{equation}\label{symeqs}
    \sum_{i=0}^{k-l}(-1)^{k-i}
    \begin{pmatrix}
    k-i \\
    l \\
    \end{pmatrix}
    \big(F_{k-i}\cdot W\big)^{(k-i-l)}=W\cdot F_l^*,\quad
    l=0,\ldots,k.
    \end{equation}
Then the differential operator $D_k$ (defined in (\ref{dop})) is
symmetric with respect to the weight matrix $W$.
\end{teor}
(In (\ref{moment}) we are using the notation for the falling or
bounded factorial $(x)_n$ defined by $(x)_n=x(x-1)\cdots(x-n+1)\quad
\mbox{for}\quad n>0,\; (x)_0=1$.)
\bigskip \begin{proof}The first and second parts are shown to be equivalent in
Proposition 4 in \cite{DdI}. Likewise, the last part can be found in
Theorem 5 in \cite{DdI}.

\end{proof}

Two advantages of the moment equations (\ref{moment}) are that
they are equivalent to the symmetry of $D_k$ and that there are no additional assumptions on the weight matrix $W$ (in
particular, neither smoothness nor boundary conditions
are required). However, even if
one can solve these moment equations it can be very difficult to
recover the weight matrix $W$ from its moments. Besides that the
moment equations turn out to become a suitable tool in all the
examples throughout this paper since those equations are going to be
the key to characterize the convex cone $\Upsilon (D)$ associated with
some of the differential operators $D$ considered in the next
section. The moment and symmetry equations for differential
operators of order two appeared for the first time in \cite{D1} when
$F_0W=WF_0^*$ and later in \cite{DG1} (see also \cite{GPT} in the case of the symmetry equations). As we remarked
in the Introduction, most of the examples which appeared in the last
years have been obtained solving these equations, while some others
came from group representation theory.

\bigskip

We now present a result that will be used to generate examples of
second order differential operators having infinitely many families
of orthogonal matrix polynomials as eigenfunctions. The idea is to
find certain constraints which guarantee the symmetry of a
differential operator with respect to both $W$ and a new weight
matrix obtained from $W$ by adding a Dirac distribution at one
point.

Let $W$ be a weight matrix and consider
\begin{equation}\label{Wd}
    \widetilde{W}(t)=W(t)+\delta_{t_0}(t)M(t_0),
\end{equation}
where $\delta_{t_0}(t)=\delta(t-t_0)$ is the Dirac delta
distribution or the ``impulse symbol'' introduced by P. A. M. Dirac
in \cite{D}, which we will consider as a measure. The Hermitian
positive semidefinite matrix $M(t_0)$ depends on the point where the
Dirac distribution is added.

Weight matrices of the form (\ref{Wd}) were considered in \cite{YMP,
YMP2} (to study asymptotic properties of the corresponding modified
Jacobi matrix) for $W$ in the Nevai class, i.e.,
with convergent recurrence coefficients.

The moments of $\widetilde W$ are related with the moments of $W$ by
the formula
\begin{equation*}\label{mom}
    \widetilde{\mu}_n=\int t^n d\widetilde{W}(t)=\mu_n+\int t^n
    \delta_{t_0}(t)M(t_0)dt=\mu_n+t_0^nM(t_0),\quad n=0,1,\ldots
\end{equation*}
Observe that in the special case of $t_0=0$ the only modified moment
is the first one $\widetilde{\mu}_0=\mu_0+M(0)$, and then
$\widetilde{\mu}_n=\mu_n$ for $n=1,2,\ldots$

The following theorem gives conditions for the symmetry of a
differential operator $D_k$ with respect to the weight matrices $W$
and $W+\delta_{t_0}M(t_0)$.

\begin{teor}\label{TNC}
Let $D_k$ be a differential operator of order $k$ as in (\ref{dop}).
Let $W$ be a weight matrix. Assume that associated with the real
point $t_0\in \RR$ there exists a Hermitian positive semidefinite
matrix $M(t_0)$ satisfying
\begin{align}\label{cond}
    F_j(t_0)M(t_0)&=0, \quad j=1,\ldots,k, \\
     \nonumber F_0M(t_0)&=M(t_0)F_0^*.
\end{align}
Then the operator $D_k$ is symmetric with respect to $W$ if and only
if it is symmetric with respect to
$\widetilde{W}=W+\delta_{t_0}M(t_0)$.
\end{teor}

\begin{proof} Recalling definitions around (\ref{moment}) for
$\widetilde{W}=W+\delta_{t_0}M(t_0)$, we produce
$$
\widetilde{B}_n^l=\sum_{i=0}^l
F_{l-i}^l\widetilde{\mu}_{n-i}=B_n^l+t_0^{n-l}F_l(t_0)M(t_0),\quad
l=0,\ldots,k.
$$
Using conditions (\ref{cond}) for $j=1,\ldots,k$, we obtain
$$
\widetilde{B}_n^0=B_n^0+t_0^{n}F_0M(t_0), \quad
\widetilde{B}_n^l=B_n^l,\quad l=1,\ldots,k.
$$
Consequently, this shows that equations (\ref{moment}), $l=1,\ldots
, k$, are just the same for $W$ and $\widetilde{W}$. For $l=0$,
equations (\ref{moment}) for $W$ and $\widetilde{W}$ are,
respectively:
\begin{align*}
\sum_{i=0}^{k-1}(n)_{k-i}B_n^{k-i}+B_n^0&=(B_n^0)^*,\\
\sum_{i=0}^{k-1}(n)_{k-i}B_n^{k-i}+B_n^0+t_0^{n}F_0M(t_0)&=(B_n^0)^*+t_0^{n}M(t_0)F_0^*.
\end{align*}
The last condition in (\ref{cond}) shows again that those equations
are the same for $W$ and $\widetilde{W}$.

\end{proof}

Note that for $N=1$ (i.e. in the scalar case) Theorem \ref{TNC} implies that either $M=0$ or there exists a common
zero for all coefficients of the differential operator. For
instance, for $k=2$ there is no such a common zero for the
classical families of Hermite, Laguerre and Jacobi.

For a weight matrix $W$ the constraints (\ref{cond}) mean that by
adding a Dirac distribution to $W$ the chances of $W$ and $W+
\delta_{t_0}M(t_0)$ sharing a symmetric differential operator of
order $k$ increase with the number of linearly independent symmetric
differential operators of order $k$ for $W$. In fact, all the
examples we show in the next section are built from a weight matrix
$W$ having several linearly independent second order differential
operators.

Let us notice that once we generate $W$, $D_k$, $t_0$ and $M(t_0)$
satisfying constraints (\ref{cond}) of the theorem above we can
produce not only one single weight matrix for which $D_k$ is
symmetric, but also a two dimensional convex cone of weight matrices
for which $D_k$ is symmetric as well. If Theorem \ref{TNC} holds for
$W+ \delta_{t_0}M(t_0)$ then automatically also holds for $\gamma
W+\zeta  \delta_{t_0}M(t_0)$ where $\gamma>0$  and $\zeta \ge 0$.
All the examples of weight matrices in $\Upsilon (D_k)$ which we
show in this paper differ (up to a multiplicative constant) in a
Dirac distribution. We think that this is by no means a general
result, and are confident that other different situations may occur
(such as the existence of differential operators $D_k$ for which
$\Upsilon (D_k)$ contains a two parametric family of weight matrices
absolutely continuous with respect to the Lebesgue measure).

We will characterize the convex cone $\Upsilon (D)$ for some of the
examples in the next section.  To do that we will need the following
result (which it is interesting in its own right). Consider the
Fourier transform $\mathcal F(W)$ of a weight matrix $W$ defined by
$$
\mathcal F(W)(x)=\int _\RR e^{itx}dW(t), \quad x\in \RR.
$$

\begin{lema}\label{nuevo} Assume that the weight matrix $W$ has moments $(\mu _n)_n$ satisfying
\begin{align}\label{cmmd}
\lim _m\tr (\mu_{2m})\frac{r^{2m}}{(2m)!}=0, \quad \mbox{for all}
\quad 0<r<R, \quad R>0,
\end{align}
where $\tr (X)$ stands for the trace of the matrix $X$. Then there
exists an analytic function $\Phi$ in the strip $\{ z\in \CC : \vert
\Im z\vert <R\}$ such that $\Phi(x)=\mathcal F(W)(x)$, $x\in \RR$.
\end{lema}

\begin{proof}

For any positive semidefinite matrix $X$ it is straightforward that
$0\preccurlyeq X\preccurlyeq \tr (X) I$ (where $\preccurlyeq$ stands
for the usual positive semidefinite ordering). Using the
Cauchy-Schwarz inequality we have that
\begin{align*}
\vert uXv^*\vert \le \vert uX^{1/2}X^{1/2}v^*\vert\le
(uXu^*)^{1/2}(vXv^*)^{1/2}\le \Vert u\Vert \Vert v\Vert \tr(X),
\end{align*}
for any vectors $u,v\in \CC ^N$ (where $\Vert u\Vert =\sqrt
{uu^*}$).

For a weight matrix $W$ ($\tr (W)$ is a positive measure), the
previous inequality gives
\begin{align}\label{traza}
\vert uWv^*\vert \le \Vert u\Vert \Vert v\Vert \tr(W),
\end{align}
for any vectors $u,v\in \CC ^N$.

The Fourier transform $\mathcal F(W)$ of a weight matrix $W$ is a
$\mathcal C^\infty$ function in $\RR$ and
\begin{align}\label{dertf}
\mathcal F(W)^{(n)}(x)=\int _\RR (it)^ne^{itx}dW(t), \quad x\in
\RR,\quad n\in\mathbb{N}.
\end{align}

Fixing a number $a\in \RR $, the Lagrange remainder for the Taylor
formula gives for each $x\in \RR$ and $k\ge 0$ a real number
$y_{k,x}$ for which
$$
\mathcal F(W)(x)=\sum _{n= 0}^{k-1}\mathcal
F(W)^{(n)}(a)\frac{(x-a)^n}{n!}+\mathcal
F(W)^{(k)}(y_{k,x})\frac{(x-a)^k}{k!}.
$$
Assume that we have proved that for all $x\in \RR $, $\vert x-a\vert
<R$,
\begin{equation}\label{seult}
\lim _k\mathcal F(W)^{(k)}(y_{k,x})\frac{(x-a)^k}{k!}=0.
\end{equation}
The power series
$$
\sum _{n= 0}^{\infty }\mathcal F(W)^{(n)}(a)\frac{(x-a)^n}{n!}
$$
defines an analytic function $\Phi_a$ in $\{ z\in \CC :\vert
z-a\vert <R\} $ satisfying $\Phi_a(x)=\mathcal F(W)(x)$, $x\in
\RR $, $\vert x-a\vert <R$. The Lemma now follows easily by using a
standard process of analytic continuation.

Let us now prove (\ref{seult}). This is equivalent to proving that for
any vectors $u,v\in \CC ^N$
\begin{equation}\label{seult2}
\lim _k u \mathcal F(W)^{(k)}(y_{k,x})v^*\frac{(x-a)^k}{k!}=0,\quad
\vert x-a\vert<R.
\end{equation}
Using (\ref{dertf}) and (\ref{traza}) we have
\begin{align*}
\vert u\mathcal F(W)^{(k)}(y_{k,x})v^*\vert &=\bigg\vert \int _\RR
(it)^ke^{ity_{k,x}}udW(t)v^*\bigg\vert
\\ &\le \int _\RR \vert t\vert ^k\vert udW(t)v^*\vert \le \Vert u\Vert \Vert v\Vert \int _\RR \vert t\vert ^k d\tr(W)(t).
\end{align*}
Hence if $k$ is even
$$
\vert u\mathcal F(W)^{(k)}(y_{k,x})v^*\vert \le \Vert u\Vert \Vert
v\Vert \tr (\mu _k);
$$
and if $k$ is odd then
$$
\vert u\mathcal F(W)^{(k)}(y_{k,x})v^*\vert \le \Vert u\Vert \Vert
v\Vert (\tr(\mu_0)+\tr(\mu _{k+1})).
$$
The limit (\ref{seult2}) follows now from (\ref{cmmd}).

\end{proof}

\bigskip
For a fixed differential operator $D$ of the form (\ref{dop1}), we
can associate another set of weight matrices
\begin{align*}\label{xd}
\X (D)=\{ W: P_n^WD=\Gamma_n P_n^W,\quad n\ge 0\},
\end{align*}
where $(P_n^W)_n$ is the sequence of  monic polynomials orthogonal
with respect to $W$.

Note that $\Upsilon (D)\subset \X (D)$ and
if $\X (D)\not =\emptyset$ then it is a cone: if $W\in \X (D)$ then
$\alpha W\in \X (D)$ for any $\alpha >0$.

In general we have $\Upsilon (D)\not =\X(D)$, as the following example shows. Let $D$ be a symmetric second order
differential operator with respect to a certain weight matrix $W$.
As a consequence the monic polynomials orthogonal with respect to
$W$ are eigenfunctions for $D$. It is clear that $iD$ is not
symmetric with respect to $W$ but  the  monic polynomials orthogonal
with respect to $W$ are still eigenfunctions for $iD$. That means
that $W\not\in \Upsilon (iD)$ but $W\in \X(iD)$, and then $\Upsilon
(iD)\not = \X(iD)$.

We are concerned that some natural questions arise regarding the
relationship between $\Upsilon (D)$ and $\X(D)$, but they are out of
the scope of this paper. We would like to quote the concluding
remark in the Introduction of \cite{CG1} because it suits very well
with this situation: ``We emphasize something that will be apparent
to any reader of this paper: the full picture of the phenomenon in
question is still far from being complete. This paper is an attempt
to describe clearly some of the new problems that one faces in the
matrix-valued case, and we pick a few examples that should give an
idea of the richness of the situation at hand''.

\section{Examples}\label{sec2}
In this section we exhibit a collection of instructive examples. The
first three examples of $(2\times2)$ weight matrices $W$ (supported
in $(-\infty,+\infty)$, $(0,+\infty)$ and $(0,1)$, respectively)
have the property that they provide four linearly independent
symmetric second order differential operators having a fixed family
of orthogonal matrix polynomials with respect to $W$ as
eigenfunctions. We show that for any real number $t_0$ we can find a
positive semidefinite matrix $M(t_0)$ and a symmetric second order
differential operator $D$ as in (\ref{difope}) satisfying the
constraints (\ref{cond}). According to Theorem \ref{TNC}, the
operator $D$ will be also symmetric with respect to $\gamma W+\zeta
\delta_{t_0}M(t_0)$, $\gamma>0,\zeta \geq0$.

The last example deals with a $(N\times N)$  weight matrix $W$
(supported in $(0,+\infty)$) with at least two linearly independent
symmetric second order differential operators. In this case we
locate $t_0$ at 0, i.e. one of the endpoints of the support of $W$,
with a mass $M$ carefully chosen.

For simplicity, our selection of examples is restricted to the field
of real numbers and to matrices with real entries, but the method we
use to find these examples is not restricted to this case.

More examples appear in the PhD dissertation of one of the authors
\cite{dI}.

\subsection{$W_a(t)=e^{-t^2}e^{At}e^{A^*t}$ with $A=\protect\begin{pmatrix}
      0 & a \\
      0 & 0
    \protect\end{pmatrix}$, $a\in \mathbb{R} \setminus \{0\}$}\label{ssH}

This weight matrix  was introduced for the first time in Section 5.1
in \cite{DG1} (for arbitrary size $N\times N$). It was deeply
explored in \cite{DG3} and a set of generators of second order
differential operators can be found in Section 6 in \cite{CG1}.

By expanding the exponential, we find that
\begin{equation}\label{peso1}
W_a(t)=e^{-t^2}e^{At}e^{A^*t}=e^{-t^2}\begin{pmatrix}
                                      1+a^2t^2 & at \\
                                      at & 1 \\
                                    \end{pmatrix}
,\quad t\in\mathbb{R},\quad a\in \mathbb{R} \setminus \{0\}.
\end{equation}

We need an expression for the
 (real) linear space of symmetric differential
operators of order at most two with respect to $W_a$. To do that we
solve equations (\ref{symeqs}) for $k=2$. For the benefit of the
reader, we recall here these equations
\begin{align}\nonumber
F_2W&=WF_2^*,\\ \label{ultult} 2(F_2W)'&=F_1W+WF_1^*,\\ \nonumber
(F_2W)''-(F_1W)'+F_0W&=WF_0^*.
\end{align}
We then get  an expression for the 5-dimensional (real) linear space
of symmetric differential operators of order at most two with
respect to $W_a$. Then, for a fixed real number $t_0$, we solve the
equations (\ref{cond}). In this case, we find that the following
differential operator
\begin{equation}\label{dopherm}
D_{a, t_0}=\partial^2F_2(t)+\partial^1F_1(t)+\partial^0F_0,
\end{equation}
where
$$
F_2(t)=\begin{pmatrix}
      -\xi_{a,t_0}^{\mp}+at_0-at & -1-(a^2t_0)t+a^2t^2 \\
      -1 & -\xi_{a,t_0}^{\mp}+at \\
    \end{pmatrix},
$$
$$
F_1(t)=\begin{pmatrix}
      -2a+2\xi_{a,t_0}^{\mp}t & -2t_0-2a\xi_{a,t_0}^{\mp}+2(2+a^2)t \\
      2t_0 & 2(\xi_{a,t_0}^{\mp}-at_0)t \\
    \end{pmatrix},
$$
$$
F_0=\begin{pmatrix}
      \xi_{a,t_0}^{\mp}+2\D\frac{t_0}{a} & 2\D\frac{2+a^2}{a^2} \\
      \D\frac{4}{a^2} & -\xi_{a,t_0}^{\mp}-2\D\frac{t_0}{a} \\
    \end{pmatrix},
$$
and the Hermitian positive semidefinite matrix $M(t_0)$
\begin{equation*}\label{MtoH}
    M(t_0)=M(a,t_0)=\begin{pmatrix}
                      (\xi_{t_0,a}^{\pm})^2 & \xi_{t_0,a}^{\pm} \\
                      \xi_{t_0,a}^{\pm} & 1\\
                    \end{pmatrix},
\end{equation*}
where
\begin{equation*}\label{xi}
\xi_{a,t_0}^{\pm}=\frac{at_0\pm\sqrt{4+a^2t_0^2}}{2},
\end{equation*}
satisfy the constraints (\ref{cond}).

This differential operator can be obtained as a linear combination
of the second order differential operators introduced in \cite{CG1}
($D_i$, $i=1,2,3,4)$, namely
\begin{equation*}
    D_{a,t_0}=\bigg(-\xi_{a,t_0}^{\mp}+\frac{2t_0}{a}\bigg)I-\xi_{a,t_0}^{\mp}D_1-\frac{4t_0}{a}D_2+\frac{4}{a^2}D_4.
\end{equation*}

Using $\xi_{a,t_0}^{+}\xi_{a,t_0}^{-}+1=0$ it is easy to verify
that the coefficients of $D_{a, t_0}$ evaluated at $t_0$ satisfy the
conditions (\ref{cond}). Since $D_{a, t_0}$ has been chosen to be
symmetric with respect to the weight matrix (\ref{peso1}), Theorem
\ref{TNC} implies that $D_{a,t_0}$ is also symmetric with respect to
any of the following weight matrices:
\begin{equation*}\label{peso1g}
    W_{a,t_0,\gamma,\zeta}(t)=\gamma W_a(t)+\zeta\delta_{t_0}(t)M(a,t_0), \quad
    \gamma>0,\zeta\geq0.
\end{equation*}

\bigskip

We now prove that our method provides all the weight matrices in the
convex cone $\Upsilon (D_{a,t_0})$ defined in (\ref{yd}), that is
$$
\Upsilon (D_{a,t_0})=\{\gamma W_a(t)+\zeta\delta_{t_0}(t)M(a,t_0);
\quad
    \gamma>0,\zeta\geq0 \}.
$$
Theorem \ref{teorsymm} implies that the sequence of moments $(\mu
_n)_n$ of each weight matrix $U\in\Upsilon (D)$  has to satisfy the
moment equations (\ref{moment}) for $k=2$ (to simplify the notation
we remove the dependence on $a$ and $t_0$). These moment equations
are
\begin{align}\label{moment21}
F_2^2\mu_n+F_1^2\mu_{n-1}+F_0^2\mu_{n-2}-\mu_n(F_2^2)^*-\mu_{n-1}(F_1^2)^*-\mu_{n-2}(F_0^2)^*=0,\quad n\ge 2;
\end{align}
\begin{align}\label{moment22}
2(n-1)(F_2^2\mu_n+F_1^2\mu_{n-1}+F_0^2\mu_{n-2})+(F_1^1\mu_n+F_0^1\mu_{n-1})+\mu_n(F_1^1)^*+\mu_{n-1}(F_0^1)^*=0,\quad n\ge 1;
\end{align}
and
\begin{align}
\label{moment23}
n(n-1)(F_2^2\mu_n+F_1^2\mu_{n-1}+F_0^2\mu_{n-2})+n(F_1^1\mu_n+F_0^1\mu_{n-1})+F_0\mu_n-\mu_n(F_0)^*&=0,\quad n\ge 0,
\end{align}
where $F_2^2, F_1^2$, $F_0^2$, $F_1^1$, $F_0^1$ and $F_0$  are,
respectively, the coefficients of the polynomials $F_2$, $F_1$ and
$F_0$. Let us recall that $D=\partial ^2F_2+\partial F_1+\partial
^0F_0$.

From equations (\ref{moment21}) and (\ref{moment22}), we get the
following expression:
\begin{align}\label{moment24}
((n-1)F^2_2+F^1_1)\mu _n +\mu_n
((n-1)F^2_2+F^1_1)^*&=((1-n)F^2_1-F^1_0)\mu_{n-1}+\mu_{n-1}((1-n)F^2_1-F^1_0)^*
\\\nonumber & \quad\quad \quad +(1-n)(F^2_0\mu_{n-2}+
\mu_{n-2}(F^2_0)^*).
\end{align}
In our example we have that
$$
(n-1)F^2_2+F^1_1=\begin{pmatrix}
    2\xi_{a,t_0}^{\mp}  & (n-1)a^2+2(2+a^2) \\
      0 & -2\xi_{a,t_0}^{\pm} \\
    \end{pmatrix}.
$$
For a fixed $n\ge 1$ the matrices $(n-1)F^2_2+F^1_1$ and
$-((n-1)F^2_2+F^1_1)^*$ do not share any eigenvalue. This implies
that equation (\ref{moment24}) defines $\mu _n$, $n\ge 1$, in
a unique way given $\mu _0$ (see \cite{G}, p. 225).

Equation (\ref{moment23}) for $n=0$ implies that $\mu _0$ has
to satisfy
$$
F_0\mu_0=\mu _0(F_0)^*.
$$
It is just a matter of computation to see that the set of solutions
$\mu_0$ of the previous equation is formed by the first moment of
the weight matrices $W_{\gamma ,\zeta}=\gamma W+\zeta
\delta_{t_0}M$, $\gamma , \zeta \in \RR $. Hence, each weight matrix
$U\in \Upsilon (D)$ has, for certain $\gamma , \zeta \in \RR$, the
same moments $(\mu_n)_n$ as $\gamma W+\zeta \delta_{t_0}M$. From the
expression (\ref{peso1}) for $W_{\gamma ,\zeta}$, we deduce that for
$m=0,1,\ldots ,$
\begin{align}\label{seseult}
 \mu_{2m}&=\gamma \begin{pmatrix}
               h_{2m}+a^2h_{2m+2} & 0 \\
               0 & h_{2m} \\
             \end{pmatrix}+\zeta t_0^{2m}M(t_0) ,\\\label{seseult2}
 \mu_{2m+1}&=\gamma\begin{pmatrix}
                      0 & ah_{2m+2} \\
                       ah_{2m+2}& 0 \\
                    \end{pmatrix} +\zeta t_0^{2m+1}M(t_0),
\end{align}
where
\begin{equation*}\label{MH}
    h_{2m}=\Gamma\bigg(m+\frac{1}{2}\bigg)=\frac{\sqrt{\pi}(2m)!}{4^m
    m!},\quad h_{2m+1}=0,
\end{equation*}
are the Hermite moments.

Since the moment $\mu_0$ has to be positive definite and $M(t_0)$ is
singular, we deduce that $\gamma >0$.

We now prove that $U=\gamma W+\zeta \delta_{t_0}M$ and
consequently $\zeta \ge 0$. This can be done using different
approaches. One of them is via Fourier transform and it works as
follows.

Equations (\ref{seseult}) and (\ref{seseult2}) show that for any $R>0$,
\begin{align*}
\lim _m\tr (\mu_{2m})\frac{R^{2m}}{(2m)!}=0.
\end{align*}
According to Lemma \ref{nuevo}, there exists an entire function
$\Phi$ such that $\Phi(x)=\mathcal F(U)(x)$, $x\in \RR$. From
(\ref{peso1}) the Fourier transform of $\mathcal F(W_{\gamma
,\zeta})$ can be computed explicitly (\cite{Le}, (4.11.4)):
$$
\mathcal F(W_{\gamma
,\zeta})(x)=\gamma\sqrt{\pi}e^{-x^2/4}\begin{pmatrix}
                                                             1+a^2(2-x^2)/4 & ix/2 \\
                                                             ix/2 & 1 \\
                                                           \end{pmatrix}+\zeta
                                                           e^{it_0x}M.
$$
This shows that $\mathcal F(W_{\gamma ,\zeta})$ is actually an
entire function.

From (\ref{dertf}) one can see that the moments $(\mu _n)_n$ of a
weight matrix $W$ are, up to a multiplicative constant, the
derivatives at $0$ of its Fourier transform: $\mu _n=\mathcal
F(W)^{(n)}(0)/i^n$, $n\ge 0$. Since $U$ and $W_{\gamma,\zeta}$ have
the same moments, we conclude that the corresponding Fourier transforms
have  at $0$ the same derivatives of any order. That is, the entire
functions $\Phi$ and $\mathcal F(W_{\gamma ,\zeta})$ are equal and
then $\mathcal F(W_{\gamma ,\zeta})(x)= \mathcal F(U)(x)$, $x\in
\RR$. So $U=W_{\gamma ,\zeta}$, in which case $\zeta$ has to be
bigger than or equal to $0$ (since $U$ is a weight matrix).

\subsection{$W_{a,\alpha}(t)=t^{\alpha}e^{-t}t^{B}t^{B^*}$ with $B=\protect\begin{pmatrix}
      1 & a \\
      0 & 0
    \protect\end{pmatrix}$, $a\in \mathbb{R} \setminus \{0\}$}\label{ssL}

This weight matrix  was introduced for the first time in Section 6.2
in \cite{DG1} (for arbitrary size $N\times N$) and it was
extensively studied in \cite{DL}. Unlike the first example, its
algebra of differential operators has not been studied in depth, but
as we mentioned at the beginning of this section, by solving
equations (\ref{ultult}), one finds that there are four linearly
independent symmetric second order differential operators.

By computing $t^B$ for  $B=\begin{pmatrix}
      1 & a \\
      0 & 0 \\
    \end{pmatrix}$ we find that
\begin{equation}\label{peso2}
W_{a,\alpha}(t)=t^{\alpha}e^{-t}t^{B}t^{B^*}=t^{\alpha}e^{-t}\begin{pmatrix}
                                                                 t^2+a^2(t-1)^2 & a(t-1) \\
                                                                 a(t-1) & 1 \\
                                                               \end{pmatrix},
\quad t\in(0,+\infty),\quad \alpha>-1.
\end{equation}

Proceeding as in Section \ref{ssH} for a fixed real number $t_0$,
we find that the following differential operator (symmetric with
respect to $W_{a,\alpha}$)
\begin{equation*}\label{doplag}
D_{a,\alpha,t_0}=\partial^2F_2(t)+\partial^1F_1(t)+\partial^0F_0,
\end{equation*}
where
\begin{align*}
\nonumber F_2(t)=&\begin{pmatrix}
                     at_0 & a^2t_0 \\
                     -t_0 & -at_0 \\
                   \end{pmatrix}
    +t \begin{pmatrix}
   \phi^{\pm}-\frac{1+(\alpha+t_0)(1+a^2)}{a} & -(\alpha+t_0+1)(1+a^2) \\
   0 & \phi^{\pm}+a \\
    \end{pmatrix}\\
    &+t^2\begin{pmatrix}
           0 & (1+a^2)(1+\alpha) \\
           0 & 0 \\
         \end{pmatrix},
\end{align*}
\begin{align*}
\nonumber F_1(t)=&\begin{pmatrix}
                     \begin{pmatrix}
                       -a(3t_0-2+(\alpha+1)^2) \\
                      +(\alpha+3)(\phi^{\pm}-\frac{t_0+\alpha+1}{a})\\
                     \end{pmatrix}
                      & \begin{pmatrix}
                          2a\phi^{\pm}-(1+a^2)(\alpha^2+2t_0+3\alpha) \\
                          +\alpha t_0(a^2-1)-(t_0+\alpha+3) \\
                        \end{pmatrix}
                       \\
                     t_0-\alpha-1 & \phi^{\pm}(\alpha+1)-a\alpha t_0 \\
                   \end{pmatrix}\\
    &+t\begin{pmatrix}
     -\phi^{\pm}+a(t_0-1)+\frac{t_0+\alpha+1}{a}& (1+a^2)(\alpha^2+3\alpha-t_0\alpha+2)+2(\alpha-1)\\
     0&-\phi^{\pm}+a\alpha\\
    \end{pmatrix},
\end{align*}
\begin{equation*}
F_0=\begin{pmatrix}
      -\frac{\phi^{\pm}}{2}+\frac{a(t_0-1)}{2}+\frac{(\alpha+2)(t_0+\alpha+1)}{2a}-\frac{a(1+\alpha)}{1+a^2}
      & 1+a\alpha\phi^{\pm}-\alpha(t_0-1)(a^2+\alpha+2)+\frac{1+\alpha}{1+a^2} \\
      \frac{1+\alpha}{1+a^2} & \frac{\phi^{\pm}}{2}-\frac{a(t_0-1)}{2}-\frac{(\alpha+2)(t_0+\alpha+1)}{2a}+\frac{a(1+\alpha)}{1+a^2} \\
    \end{pmatrix},
\end{equation*}
and the Hermitian positive semidefinite matrix $M(t_0)$
\begin{equation*}\label{Mt0L}
    M(t_0)=M(a,\alpha,t_0)=\begin{pmatrix}
                      (\phi_{a,\alpha,t_0}^{\pm})^2 & \phi_{a,\alpha,t_0}^{\pm} \\
                      \phi_{a,\alpha,t_0}^{\pm} & 1 \\
                    \end{pmatrix},
\end{equation*}
where
\begin{equation*}\label{zeta}
\phi^{\pm}=\phi_{a,\alpha,t_0}^{\pm}=\frac{1}{2}\frac{(a^2+1)(t_0+\alpha)-a^2+1\pm\sqrt{(a^2+1)(a^2(t_0-\alpha-1)^2+(t_0+\alpha+1)^2)}}{a},
\end{equation*}
satisfy the constraints (\ref{cond}).

Thus, Theorem \ref{TNC} implies that $D_{a,\alpha,t_0}$ is symmetric
with respect to any of the following weight matrices:
\begin{equation*}\label{peso1lag}
    W_{a,\alpha,t_0,\gamma, \zeta}(t)=\gamma W_{a,\alpha}(t)\chi_{(0,+\infty)}(t)
    + \zeta \delta_{t_0}(t)M(a,\alpha,t_0), \quad t\in\mathbb{R},\quad \alpha>-1,\quad \gamma>0,\zeta\geq0.
\end{equation*}
Note that the discrete mass on the Delta distribution can be located
in or out of the support of the original weight matrix
(\ref{peso2}).

Our method provides all the weight matrices in the convex
cone $\Upsilon (D_{a,\alpha, t_0})$, i.e.
$$
\Upsilon (D_{a,\alpha,t_0})=\{W_{a,\alpha,t_0,\gamma, \zeta}(t);
\quad
    \gamma>0,\zeta\geq0 \}.
$$
This can be proved in a similar way as for the previous example. We
have that
$$
(n-1)F^2_2+F^1_1=\begin{pmatrix}
    \phi^{\mp}-a\alpha & (n-1)(1+\alpha)(1+a^2) \\
      0 &  -\phi^{\pm}+a\alpha \\
    \end{pmatrix},
$$
hence, for a fixed $n\ge 1$ again the matrices $ (n-1)F^2_2+F^1_1$
and $-((n-1)F^2_2+F^1_1)^*$ do not share any eigenvalue.

From the expression (\ref{peso2}) for $W_{\gamma ,\zeta}$ we deduce
that for $n=0,1,\ldots$
\begin{equation}\label{ML}
    \mu_{n}=\gamma \Gamma(n+\alpha+1)\begin{pmatrix}
                \theta_n& a(\alpha+n) \\
                a(\alpha+n) & 1 \\
             \end{pmatrix}+\zeta t_0^nM(t_0)
\end{equation}
where
$$
\theta_n=a^2(\alpha^2+(2n+1)\alpha+n^2+n+1)+\alpha^2+(2n+3)\alpha+(n+1)(n+2),
$$
and that $\mathcal F(W_{a,\alpha,t_0,\gamma, \zeta})$ is analytic
 in the cut plane $\CC \setminus \{ -ix:x\in \RR , 1<x \}$ (the explicit expression
of $\mathcal F(W_{a,\alpha,t_0,\gamma, \zeta})$ can be computed
easily from \cite{Le}, (1.5.1)).

Equation (\ref{ML}) shows that for all $0<r<1$
\begin{align*}
\lim _{m}\tr (\mu_{2m})\frac{r^{2m}}{(2m)!}=0.
\end{align*}
And we can proceed as in the previous example applying Lemma
\ref{nuevo} in the strip $\{ z\in \CC : \vert \Im z \vert <1\}$.

\subsection{An example supported in $(0,1)$}\label{ssJ}

The following weight matrix is a modification of the one introduced
in \cite{PT1}:
\begin{equation*}\label{peso3}
    W_{\alpha,\beta,k}(t)=t^{\alpha}(1-t)^{\beta}\begin{pmatrix}
                                           kt^2+\beta-k+1 & (\beta-k+1)(1-t) \\
                                           (\beta-k+1)(1-t) & (\beta-k+1)(1-t)^2 \\
                                         \end{pmatrix},\quad
                                         t\in(0,1),
\end{equation*}
where $\alpha,\beta>-1$ and $0<k<\beta+1$.

The difference between $W_{\alpha,\beta,k}$ and the weight matrix
introduced in \cite{PT1} ($N=2$), is that $W_{\alpha,\beta,k}$
enjoys four linearly independent symmetric second order differential
operators while that in \cite{PT1} enjoys only two. The example in
\cite{PT1} is related to a group-theoretical situation but this is
not the case of our $W_{\alpha,\beta,k}$, as far as we know.

As in the previous examples, for a fixed real number $t_0$ (except
for $t_0=-\frac{1}{k}(\alpha+\beta-k+2)$), no matter if it is
located in or out of the support of the weight matrix $W_{\alpha,
\beta, k}$, we find a symmetric second order differential operator
$D_{\alpha,\beta ,k,t_0}$ with respect to $W_{\alpha, \beta , k}$
and a Hermitian positive semidefinite matrix $M(t_0)$ satisfying the
constrains (\ref{cond}). Hence, this operator $D_{\alpha,\beta
,k,t_0}$ is also symmetric for any of the weight matrices $\gamma
W_{\alpha, \beta, k}\chi_{(0,1)}+\zeta \delta_{t_0}M(t_0)$,
$\gamma>0, \zeta \geq0$. The convex cone $\Upsilon (D_{a,\alpha,
\beta, t_0})$ is formed by those weight matrices. This can be proved
as in the previous examples. Since in this example the weight
matrices have compact support $[0,1]\cup \{ t_0\}$, it
follows that the moments $(\mu _n)_n$ satisfy
$\tr(\mu_n)\le \gamma \tr (\mu _0)+\tr (M(t_0))\vert\zeta \vert\vert
t_0\vert ^n$, $n\ge 0$. Then by Lemma \ref{nuevo} the Fourier
transforms are entire functions.

Since the formulas for arbitrary $\alpha ,\beta$ and $k$ are very
long we show here only one concrete example: $\alpha=0$, $\beta =0$
and $k=1/2$:
\begin{align*}
F_2(t)&= \begin{pmatrix}
         (1-t)(-t_0+t(2t_0-3))+\D\frac{t(1-t)(1-t_0)}{\varphi^{\pm}} & 2t+t_0-2t_0t-t_0t^2 \\
         -t_0(1-t)^2 & (t-1)(t-t_0)+\D\frac{t(1-t)(1-t_0)}{\varphi^{\pm}} \\
       \end{pmatrix}, \\
F_1(t)&= \begin{pmatrix}
         12t-10+8t_0(1-t)+\D\frac{(t_0-1)(4t-3)}{\varphi^{\pm}} & 9-t-4t_0(t+2)+\D\frac{2(t_0-1)}{\varphi^{\pm}} \\
         (1-4t_0)(t-1) & 4(t-t_0)+\D\frac{(t_0-1)(4t-1)}{\varphi^{\pm}} \\
       \end{pmatrix},\\
F_0&= \begin{pmatrix}
         \D\frac{(1-t_0)(1+2\varphi^{\pm})}{2\varphi^{\pm}} & -t_0-3+\D\frac{1-t_0}{2\varphi^{\pm}} \\
         -t_0+\D\frac{1-t_0}{2\varphi^{\pm}} & -\D\frac{(1-t_0)(1+2\varphi^{\pm})}{2\varphi^{\pm}} \\
       \end{pmatrix},
\end{align*}
and
\begin{equation*}\label{Mt0J}
    M(t_0)=\begin{pmatrix}
                      1 & \varphi^{\pm} \\
                      \varphi^{\pm} & (\varphi^{\pm})^2 \\
                    \end{pmatrix},
\end{equation*}
where
\begin{equation*}
\varphi^{\pm}(t_0)=\frac{2-t_0\pm\sqrt{2t_0^2-2t_0+1}}{t_0+3}.
\end{equation*}
We need to impose $t_0\neq -3$ to avoid singularities in $M(t_0)$
($t_0=1$ gives $\varphi^{-}(1)=0$, in which case
$(1-t_0)/\varphi^{-}$ has to be taken equal to $1$ in the entries of
$F_2, F_1$ and $F_0$).

\subsection{An example of arbitrary size}\label{Lagnew}

We consider here a weight matrix defined by
\begin{equation}\label{pesolagnew}
W(t)=t^{\alpha}e^{-t}e^{At}t^{\frac{1}{2}J}t^{\frac{1}{2}J^*}e^{A^*
t},\quad t\in(0,\infty),\quad \alpha>-1,
\end{equation}
where
\begin{equation}\label{AAA}
J=\sum_{i=1}^{N} (N-i)E_{ii},\quad  A=\sum_{i=1}^{N-1} \nu_i
E_{i,i+1},\quad \nu_i\in\mathbb{R}\setminus\{0\}.
    \quad i=1,\ldots,N-1.
\end{equation}
Here we are using $E_{ij}$ to denote the matrix with entry $(i,j)$
equal to 1 and 0 otherwise.

This weight matrix  was introduced for the first time in \cite{DdI}.
This example enjoys the special property of having symmetric odd
order differential operators, a phenomenon that is not possible in
the classical scalar theory. For more details, the reader should
consult \cite{DdI}.

It is proved in \cite{DdI} that (\ref{pesolagnew}) always has a
symmetric second order differential operator given by
\begin{equation*}\label{l21lag}
    D_1=\partial^2tI+\partial^1[(\alpha+1)I+J+t(A-I)]+\partial^0[(J+\alpha I)A-J],
\end{equation*}
where $A$ and $J$ are defined in (\ref{AAA}). Note that there are
$N-1$ free parameters in $D_1$. Assuming the following conditions on
the parameters $\nu_1,\ldots,\nu_{N-2}$:
\begin{equation*}
i(N-i)\nu_{N-1}^2=(N-1)\nu_{i}^2+(N-i-1)\nu_{i}^2\nu_{N-1}^2, \quad
i=1,\ldots,N-2,
\end{equation*}
the weight matrix $W$ has another symmetric second order
differential operator
\begin{equation*}\label{l22lag}
    D_2=\partial^2G_2(t)+\partial^1G_1(t)+\partial^0G_0,
\end{equation*}
where the coefficients are given by
\begin{align*}
  G_2(t)&=t(J-At), \\
  G_1(t)&=((1+\alpha)I+J)J+Y-t(J+(\alpha+2)A+Y^*-AY+YA),\\
  G_0&=\D\frac{N-1}{\nu_{N-1}^2}[J-(\alpha I+J)A],
\end{align*}
$A$ and $J$ are defined in (\ref{AAA}), and $Y=\D\sum_{i=1}^{N-1}
\frac{i(N-i)}{\nu_i}E_{i+1,i}$. Note that now the only free
parameter is $\nu_{N-1}$.

Let us define the Hermitian positive semidefinite matrix (and
singular) $M=v^*v$ , where $v$ is a row vector
\begin{equation*}
    v=\D\sum_{j=1}^{N-1}
    \bigg(\prod_{k=1}^{N-j}\D\frac{\nu_{N-k}(\alpha+k)}{k}\bigg)e_{j}+e_N,
\end{equation*}
where $e_j=(0,\ldots,0,1,0,\ldots,0)$ denotes the $j$-th unit
vector. A simple computation gives
\begin{equation*}\label{Mt0Lnew}
    M=\D\sum_{i,j=1}^N \bigg(\prod_{k=\min\{i,j\}}^{\max\{i,j\}-1}\D\frac{\nu_k(\alpha+N-k)}{N-k}\bigg)
    \bigg(\prod_{k=1}^{N-\max\{i,j\}}\D\frac{\nu_{N-k}(\alpha+k)}{k}\bigg)^2E_{ij}
\end{equation*}
(where for $m>n$ we take $\prod_{k=m}^n=1$).

Let $D$ be the following differential operator $D=-(N-1)D_1+D_2.$ As
a linear combination of symmetric operators with respect to the
weight matrix (\ref{pesolagnew}), $D$ is also symmetric with respect
to this weight matrix. The differential coefficients ${F}_i(t),\;
i=0,1,2,$ of $D$ are
\begin{align*}
  F_2(t)&=t(-(N-1)I+J-At), \\
  F_1(t)&=Y-\sum_{i=1}^N (i-1)(\alpha+N-i+1)E_{ii}+t(J-(\alpha+N+1)A-Y^*),\\
  F_0&=\D\frac{(N-1)(1+\nu_{N-1}^2)}{\nu_{N-1}^2}[J-(\alpha I+J)A].
\end{align*}
Evaluating them at $t=0$ and considering the bidiagonal structure
of $F_1(0)$ and $F_0$, it is easy to check that $F_1(0)v^*=0$ and
$F_0v^*=0$. Thus the conditions (\ref{cond}) are satisfied
($F_2(0)=0$).  Hence Theorem \ref{TNC} implies that $D$ is
symmetric with respect to any of the following weight matrices:
\begin{equation*}\label{peso1lagnew}
   W_{\gamma,\zeta}(t)=\gamma W(t)+ \zeta \delta(t)M,
\quad t\in(0,+\infty),\quad \alpha>-1,\quad\gamma>0, \zeta \geq0.
\end{equation*}

We illustrate the case of size $2\times2$ ($\nu_1=a$). The weight
matrix $W_{\gamma,\zeta}$ is
\begin{equation*}
    W_{\gamma, \zeta }(t)=\gamma t^{\alpha}e^{-t}\begin{pmatrix}
                                      t(1+a^2t) & at \\
                                      at & 1 \\
                                    \end{pmatrix}+\zeta\begin{pmatrix}
                      a^2(\alpha+1)^2 & a(\alpha+1) \\
                      a(\alpha+1) & 1 \\
                    \end{pmatrix}\delta(t),
\end{equation*}
while the second order differential operator $D$ is
\begin{align*}
D=\partial^2\begin{pmatrix}
      0 & -at^2 \\
      0 & -t \\
    \end{pmatrix}
    +\partial^1\begin{pmatrix}
      t & -\D\frac{(1+a^2(\alpha+3))t}{a} \\
      \D\frac{1}{a} & -(\alpha+1) \\
    \end{pmatrix}
    +\partial^0\begin{pmatrix}
      \D\frac{a^2+1}{a^2}& -\D\frac{(1+a^2)(\alpha+1)}{a} \\
      0 & 0 \\
    \end{pmatrix}.
\end{align*}

\bigskip
\textit{Acknowledgements:} The authors would like to thank the
anonymous referee for the careful reading of this paper as well as
comments and suggestions.

\end{document}